\documentclass{article}
\usepackage{graphicx}

\usepackage{amsmath}
\usepackage{url}
\usepackage[T2A]{fontenc}
\usepackage[russian,english,british]{babel}

\begin{document}

\title{Inevitable Dottie Number. Iterals of cosine and sine}
\author{Valerii Salov}
\date{}
\maketitle

\begin{abstract}
The unique real root of $\cos(x) = x$, recently referred to as the Dottie number, is expressed as an iteral of cosine. Using the derivatives of iterals, it is shown why this number is achieved starting from any real number, when the iterates of cosine successfully approach infinity, and how this affects the Maclaurin series of the iterals. Properties of the iterals of cosine and sine and their derivatives are considered. A C++ template for iteral is applied for computation of Julia sets.
\end{abstract}

\section{Introduction}

If one sets a calculator to work with radians, types any number, and presses the button $\cos$ over and over again, then soon the screen shows unchanging further decimal digits $0.739085133215\ldots$. This is a truncated root of the equation $\cos(x) = x$. Remembering such a story occurred with Dottie, a professor of French, Samuel Kaplan refers to the root as the \textit{Dottie number} \cite{kaplan}. He formulates exercises  involving the number for Advanced Calculus, Differential Equations, Chaos, Complex Variables, and concludes

\textit{It is unlikely that the Dottie number will enter the annals of great constants along side $e$, $\pi$, the Golden Ratio and many others. However, the Dottie number and its story might make good teaching elements for others out there. I also imagine there are many other interesting facets of the Dottie number yet to be discovered.}

The two properties that this is 1) the only real root of $\cos(x) = x$, and 2) the attractor in the iterates of $\cos$, where the domain of attraction is the entire real line had been known prior to the Kaplan's naming proposal.

The first property implies computation of the root. In 19th century, T. Hugh Miller \cite{miller} evaluates eight decimal places by a method of trial. He mentions that the problems of finding the angles satisfying the equation occur in analytical geometry and presents a proof of the uniqueness of the real root

\textit{For while $x$ increases from $0$ to $\frac{\pi}{2}$, $\cos(x)$ diminishes from $1$ to $0$. Also $\cos(x)$ is positive while $x$ lies between $0$ and $-1$.}

The chart of the line $y = x$ and curve $y = \cos(x)$ intersecting only in one point is a geometrical proof. Then, for this root $\cos(x) = x = \arccos(x)$. The first hundred decimal digits of the Dottie number are found in \cite{kaplan} and more in \cite{broukhis}. Hrant Arakelian presents 6,400,000 of them \cite{arakelian}. Miller proposes a method for computation of complex roots \cite{miller} and reminds that their number is infinite \cite{miller2}. In \cite{miller} he notices without a proof

\textit{It is easy to show that one value of the angle is $\cos\cos\cos\dots\cos 1$, when the cosine is taken successively to infinity}.

He does not discuss iterative properties of cosine. It is not clear, whether Miller knew that "cosine taken successfully to infinity" gives the same result for any real number. The \textit{Miller number} is a suitable name, if the community wants to emphasize his computational contribution.

The second property is discussed, for instance, by Robert Devaney \cite[pp. 39 - 56]{devaney}. Applying graphical analysis, he illustrates dynamical properties of the functions $\sqrt x, x^2, 4x(1 - x), x^2 - 2, \sin(x),$ and $\cos(x)$. Graphical analysis shows that iterates of cosine starting from arbitrary real values tend to the real root of $\cos(x) = x$. This method is used by researchers to study function iterates. For instance, Mitchel Feigenbaum discovers a number named today after him \cite{feigenbaum}.

Probably, generations using calculators have faced with the "remarkable $\cos$ experiment". Psychology of "unexpected" is responsible for the first impression. Many would agree without a calculator, that successive iterations of $\sqrt x$ tends to one beginning from any $x > 1$ and, after more thinking, that this is true for $x > 0$. But for $\cos$ it is not obvious. One is easier to memorize than $0.739085\ldots$.

Since 1980th, Hrant Arakelian emphasizes the fundamental meaning of the real root of $\cos(x) = x$ viewing it as \textit{a hidden link of mathematics} and suggesting the names \textit{cosine constant, cosine superposition constant, cosine attractor, cosine fixed point} \cite{arakelian}. He proposes to denote it by the first lower case character of the Armenian alphabet pronounced [a:] and claims relationships between this and some fundamental physical constants. If his findings will be eventually recognized, then the name \textit{Arakelian number} will become inevitable. The author is not yet ready to share Arakelian's views.

It is natural for the author \cite{salov} to apply to the iterates of $\cos(x)$ the notation for the iteration of functions recently introduced and referred by him to as the \textit{iteral}
\begin{displaymath}
\textrm{\LARGE{\CYRI}}_{v}^{n}\cos(x).
\end{displaymath}
For a function of one variable, $v$ is the initial value of the variable, and $n$ is the number of the iterations. On each iteration a previous function value becomes its new argument. By definition, the iterals of $\cos(x)$ are
\begin{displaymath}
\textrm{\LARGE{\CYRI}}_{v}^{0}\cos(x)=v, \; \textrm{\LARGE{\CYRI}}_{v}^{1}\cos(x) = \cos(v), \; \textrm{\LARGE{\CYRI}}_{v}^{2}\cos(x) = \cos(\cos(v)).
\end{displaymath}
The "intriguing calculator experiment" hints that
\begin{displaymath}
\forall v \in \mathrm{R} \; \textrm{\LARGE{\CYRI}}_{v}^{\infty}\cos(x) = \textrm{Dottie number} = D.
\end{displaymath}
Here, $\forall$ "for all" is the universal quantification of proposition, and \textrm{R} denotes the set of real numbers. The goal of this article is to review the properties of the iterals of cosine and sine, derivatives of the iterals, Maclaurin series, and show why the Dottie number is achieved.

\section{Iterals of Cosine and Sine and First Derivatives}

By the definitions of the Dottie number and iteral
\begin{displaymath}
\forall n \in \mathrm{N_0} \; \textrm{\LARGE{\CYRI}}_{D}^{n}\cos(x) = D.
\end{displaymath}
$\mathrm{N_0}$ is the set of natural numbers and zero. This reflects that $D$ is a fixed point of cosine. Since cosine is even function, the following bounds apply
\begin{displaymath}
\forall v \in \textrm{R} \; \textrm{\LARGE{\CYRI}}_{v}^{1}\cos(x) \in [-1, 1] \; \textrm{and} \; \textrm{\LARGE{\CYRI}}_{v}^{2}\cos(x) \in [\cos(1), 1].
\end{displaymath}
If the function values belong to its domain up to $n$ iterations, then \cite{salov}
\begin{displaymath}
\textrm{\LARGE{\CYRI}}_{v}^{n+1}f(x) = \textrm{\LARGE{\CYRI}}_{\textrm{\LARGE{\CYRI}}_{v}^{n}f(x)}^{1}f(x).
\end{displaymath}
For cosine and sine this is the case for all $v \in \mathrm{C}$ and $n \ge 0$, where $\mathrm{C}$ denotes the set of complex numbers. Iterals of cosine and sine are the functions of the initial value $v$. They create two families of functions, where each member corresponds to the order of iteral $n$. When the initial value $v$ is a valid expression containing the variable $x$, parameters, and functions the iteral creates new functions of $x$. The function under the iteral sign can be a vector function \cite{salov}.

Replacing $v$ with $x$ makes the iterals of cosine and sine functions of $x$
\begin{displaymath}
\textrm{\LARGE{\CYRI}}_{x}^{n}\cos(x) \; \textrm{and} \; \textrm{\LARGE{\CYRI}}_{x}^{n}\sin(x),
\end{displaymath}
For basic theorems on iterated functions see \cite{bennett}. Let us find their first derivatives using the chain rule and iteral notation. We get for the orders $n = 0, 1, 2$
\begin{displaymath}
\begin{array}{ccccc}
n & \textrm{\LARGE{\CYRI}}_{x}^{n}\cos(x) & \left(\textrm{\LARGE{\CYRI}}_{x}^{n}\cos(x)\right)' & \textrm{\LARGE{\CYRI}}_{x}^{n}\sin(x) & \left(\textrm{\LARGE{\CYRI}}_{x}^{n}\sin(x)\right)' \\
0 & x & 1 & x & 1 \\
1 & \cos(x) & -\sin(x) & \sin(x) & \cos(x) \\
2 & \cos(\cos(x)) & \sin(x)\sin(\cos(x)) & \sin(\sin(x)) & \cos(x)\cos(\sin(x)) \\
\end{array}
\end{displaymath}
With the increasing order of the iteral of cosine the \textit{previous} first derivative multiplied by $-\sin(\textrm{\LARGE{\CYRI}}_{x}^{n-1}\cos(x))$ yields the \textit{current} first derivative. For the first derivative of the iteral of sine the multiplier is $\cos(\textrm{\LARGE{\CYRI}}_{x}^{n-1}\sin(x))$. Thus,
\begin{displaymath}
\frac{d\textrm{\LARGE{\CYRI}}_{x}^{n}\cos(x)}{dx} = (-1)^n \prod_{k=0}^{n-1} \sin(\textrm{\LARGE{\CYRI}}_{x}^{k}\cos(x)) \; \textrm{and}
\end{displaymath}
\begin{displaymath}
\frac{d\textrm{\LARGE{\CYRI}}_{x}^{n}\sin(x)}{dx} = \prod_{k=0}^{n-1} \cos(\textrm{\LARGE{\CYRI}}_{x}^{k}\sin(x)).
\end{displaymath}
For $n=0$ no multiplications are performed. The default value of $\prod$ is one. The iterals make the expressions compact for large $n$. Let us prove the formula for the iterals of cosine using mathematical induction. The proof for the iterals of sine is similar and omitted. The formula is correct for the orders 0, 1, 2. Let it be valid for $n > 2$. Then, for $n+1$ we get
\begin{displaymath}
\frac{d\textrm{\LARGE{\CYRI}}_{x}^{n+1}\cos(x)}{dx} = \frac{d\textrm{\LARGE{\CYRI}}_{\textrm{\LARGE{\CYRI}}_{x}^{n}\cos(x)}^{1}\cos(x)}{dx} = \frac{d\cos(\textrm{\LARGE{\CYRI}}_{x}^{n}\cos(x))}{dx}=
\end{displaymath}
\begin{displaymath}
=-\sin(\textrm{\LARGE{\CYRI}}_{x}^{n}\cos(x))\frac{d\textrm{\LARGE{\CYRI}}_{x}^{n}\cos(x)}{dx} =
\end{displaymath}
\begin{displaymath}
= -\sin(\textrm{\LARGE{\CYRI}}_{x}^{n}\cos(x)) (-1)^n \prod_{k=0}^{n-1} \sin(\textrm{\LARGE{\CYRI}}_{x}^{k}\cos(x)) =  (-1)^{n+1} \prod_{k=0}^{n} \sin(\textrm{\LARGE{\CYRI}}_{x}^{k}\cos(x)) =
\end{displaymath}
\begin{displaymath}
=(-1)^{h} \prod_{k=0}^{h-1} \sin(\textrm{\LARGE{\CYRI}}_{x}^{k}\cos(x)).
\end{displaymath}
This is exactly the same formula for $n$ replaced with $h=n+1$.

Since for $x \in [-1, 1]$ sine is a monotonously increasing function and $-1 < -0.84... = \sin(-1)  \leq \sin(x) \leq \sin(1) = 0.84... < 1$, each factor, with the exception of the first one, in the first derivative of a cosine iteral of the order $n > 1$ is by absolute value less than 0.85. This is because of the mentioned bounds of the cosine iterals serving as arguments for the sines. For any $x$ the iterals of the order $n$ form a variant, where each member starting from the second by absolute value is less than $|q=0.9|^n$. The latter is infinitesimal variant for $n \rightarrow \infty$ \cite[p. 50]{fihtengoltz}.  Then,
\begin{displaymath}
\forall x \in \textrm{R} \lim_{n \rightarrow \infty} \frac{d\textrm{\LARGE{\CYRI}}_{x}^{n}\cos(x)}{dx} = \lim_{n \rightarrow \infty} (-1)^n \prod_{k=0}^{n-1} \sin(\textrm{\LARGE{\CYRI}}_{x}^{k}\cos(x)) = 0.
\end{displaymath}
Because all iterals of cosine also intersect at the fixed point $D$ this ensures that
\begin{displaymath}
\forall x \in R \lim_{n \rightarrow \infty}\textrm{\LARGE{\CYRI}}_{x}^{n}\cos(x) = D.
\end{displaymath}
The idea to consider the limit of the first bound derivative of the iterals suggested by Alexander Tumanov in a private discussion simplifies the proof. The latter formulated by the author with a variant of iterals complements the Graphical Analysis proof and Miller's notice about $\cos\cos\cos\dots\cos 1$ iterated to infinity.

For $n > 0$ iterals of cosine are even functions, below \textrm{N} is the set of natural numbers without zero,
\begin{displaymath}
\forall n \in \mathrm{N} \; \textrm{\LARGE{\CYRI}}_{-x}^{n}\cos(x) = \textrm{\LARGE{\CYRI}}_{x}^{n}\cos(x).
\end{displaymath}
The iteral of the order zero $x$ is odd function. All iterals of sine are odd functions
\begin{displaymath}
\forall n \in \mathrm{N_0} \; \textrm{\LARGE{\CYRI}}_{-x}^{n}\sin(x) = -\textrm{\LARGE{\CYRI}}_{x}^{n}\sin(x).
\end{displaymath}

For $n > 0$ the iterals of cosine and sine are periodic functions. This follows from a theorem on iterates of periodic functions \cite{bennett}.  For $n = 1$ we deals with $\cos(x)$ and $\sin(x)$ and their period $2\pi$. For $n \ge 2$ the period of iterals of cosine is $\pi$. This is because $\cos(\cos(x + \pi)) = \cos(-\cos(x)) = \cos(\cos(x))$. Further composition by $\cos$ cannot change the result obtained after the two initial iterations. Increasing the order of iterals of sine does not reduce the period for $n \ge 2$ because $\sin$ is odd function. It remains equal to $2\pi$.

Independently on $n$, the real domain of each iteral of cosine and sine is $(-\infty, \infty)$. The ranges of the iterals for $n = 0$ are also $(-\infty, \infty)$. The ranges of the iterals for $n = 1$ are $[-1, 1]$. Again, this is true for real functions of cosine and sine. For $n = 2$ the ranges of the iteral of cosine and sine are $[\cos(1), 1]$ and $[-\sin(1), \sin(1)]$. Till now, the ranges are narrowing with the increasing order $n$ for both families. This remains true for $n > 2$. Indeed, on the interval $[-1, 1]$ $\sin$ is monotonously increasing function. It takes its minimum value at $-1$ and maximum value at $1$. The minimum is greater than $-1$ and the maximum is less than $1$. They form a narrower input range before the third iteration of $\sin$. Since $x < \sin(x)$ for $x \in [-1, 0)$ and $\sin(x) < x$ for $x \in (0, 1]$ a current iteral range is narrower than a previous one. Since $\sin$, being odd, is antisymmetric with respect to the origin of the coordinates, the bounds are moving towards zero so that it remains the middle of the iteral range. Thus, $\forall x \in \textrm{R} \vee n \in \textrm{N}$
\begin{displaymath}
\textrm{\LARGE{\CYRI}}_{-1}^{n-1}\sin(x) < \textrm{\LARGE{\CYRI}}_{-1}^{n}\sin(x) \le \textrm{\LARGE{\CYRI}}_{x}^{n}\sin(x)  \le \textrm{\LARGE{\CYRI}}_{1}^{n}\sin(x) < \textrm{\LARGE{\CYRI}}_{1}^{n-1}\sin(x).
\end{displaymath}
The $n$ order iterals of sine of 1 (-1) form a decreasing (increasing) variant with positive $p_n$ (negative $q_n$) members. Therefore, the variant $p_n$ ($q_n$) has a limit: non-negative $a$ (non-positive $b$). Let us consider $\sin(p_n)$ ($\sin(q_n)$). On the one hand, its value at the limit is $\sin(a)$ ($\sin(b)$). On the other hand, it is also an interal of the order $n+1$. Thus, $\sin(a) = a$ ($\sin(b) = b$). This is possible only, if $a=b=0$. The final part of the proof with "extra" sine applied to the variant is given by Alexander Tumanov in a private discussion. Hence,
\begin{displaymath}
\forall x \in R \lim_{n \rightarrow \infty}\textrm{\LARGE{\CYRI}}_{x}^{n}\sin(x) = 0.
\end{displaymath}

For the iterals of cosine the picture is richer. The range of the second order iteral is [cos(1), 1]. This becomes the input interval for the last $\cos$ evaluated for the third order iteral. On this interval $\cos$ is monotonously decreasing and positive function. Accordingly, it takes its minimum $\cos(1)$ and maximum $\cos(\cos(1))$ values at $1$ and $\cos(1)$. The range of the third order iteral is $[\cos(1), \cos(\cos(1))]$. Let us notice that this range is narrower than for the second order iteral because $\cos(\cos(1)) < 1$ and the low bounds for the second and third orders are identical. By these reasons, $\cos$ on this interval is also monotonously decreasing and positive. Meanwhile, this range becomes the input interval for the last $\cos$ evaluated for the fourth order iteral. It takes its minimum $\cos(\cos(\cos(1)))$ at the most right argument $\cos(\cos(1))$ and its maximum $\cos\cos(1))$ at the left most argument $\cos(1)$. The range of the fourth order iteral is $[\cos(\cos(\cos(1))),\cos\cos(1))]$. Let us notice that this range is narrower than for the third order iteral because $\cos(1) < \cos(\cos(\cos(1)))$ and, now, the top bounds for the third and fourth orders are identical. Then, on this input interval for the last $\cos$ of the fifth order iteral $\cos$ remains monotonously decreasing and positive. During two consecutive steps either low or bottom bounds are identical and then move towards each other. On any closed subinterval of $[\cos(1), 1]$ $\cos$ evaluated at the most right end is less than $\cos$ evaluated at the most left end. A range of the order $n > 1$ iteral of cosine can be expressed as
\begin{displaymath}
\left[\textrm{\LARGE{\CYRI}}_{1}^{n-1-(n \% 2)}\cos(x), \; \textrm{\LARGE{\CYRI}}_{1}^{n-2+(n \% 2)}\cos(x)\right],
\end{displaymath}
where $\%2$ is division modulo $2$. The three sequences clarify the formulas
\begin{displaymath}
\begin{array}{rrrrrrrrrr}
n & 2 & 3 & 4 & 5 & 6 & 7 & 8 & 9 & \ldots\\
n-1-(n\%2) & 1 & 1 & 3 & 3 & 5 & 5 & 7 & 7 & \ldots\\
n-2+(n\%2) & 0 & 2 & 2 & 4 & 4 & 6 & 6 & 8 & \ldots
\end{array}
\end{displaymath}
yielding
\begin{displaymath}
\forall x \in \textrm{R} \vee n > 1 \; \textrm{\LARGE{\CYRI}}_{1}^{n-1-(n \% 2)}\cos(x) \le \textrm{\LARGE{\CYRI}}_{x}^{n}\cos(x) \le \textrm{\LARGE{\CYRI}}_{1}^{n-2+(n \% 2)}\cos(x).
\end{displaymath}
Being iterals of cosine, the bounds tend to the Dottie number with $n \rightarrow \infty$.

The locations of minimums and maximums of the iterals of cosine are the roots of
\begin{displaymath}
\frac{d\textrm{\LARGE{\CYRI}}_{x}^{n}\cos(x)}{dx} = (-1)^n \prod_{k=0}^{n-1} \sin(\textrm{\LARGE{\CYRI}}_{x}^{k}\cos(x)) = 0.
\end{displaymath}
From all $\sin$ factors in the product only $\sin(x)$ and $\sin(\cos(x))$ can get zero values. For $n > 1$ the iterals of cosine in the remaining $\sin$ factors have the ranges within the interval $[\cos(1), 1]$. On this interval $\sin(x) > 0$. Thus, the minimums and maximums of the $n > 1$ iterals of cosine at the roots of $\sin(x)\sin(\cos(x)) = 0$ are equal to $\pm \frac{k\pi}{2}, \; k \in \mathrm{N_0}$.

The locations of minimums and maximums of the iterals of sine are the roots of
\begin{displaymath}
\frac{d\textrm{\LARGE{\CYRI}}_{x}^{n}\sin(x)}{dx} = \prod_{k=0}^{n-1} \cos(\textrm{\LARGE{\CYRI}}_{x}^{k}\sin(x))= 0.
\end{displaymath}
For the iterals of sine of the order $n > 0$ each expression gets the factor $\cos(x)$. Therefore, the locations coincide with those of sine because these iterals have the same period as sine. They are equal to $\frac{\pi}{2} \pm k\pi, k \in N_0$

The following figures depict the first seven iterals of cosine and sine, locations and sizes of maximums and minimums and intersection points. The unbound iteral of the order zero is the straight line y = x for both cosine and sine. It is omitted in order to maintain suitable drawing scale. Due to evenness and periodicity of the iterals of cosine, all intersection points correspond to $y = D$.

\includegraphics[width=130mm]{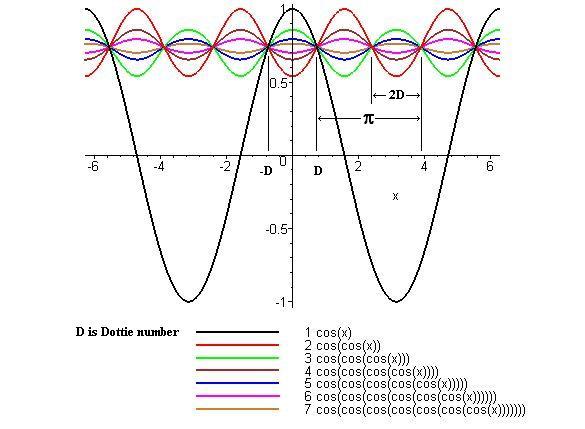}

\includegraphics[width=130mm]{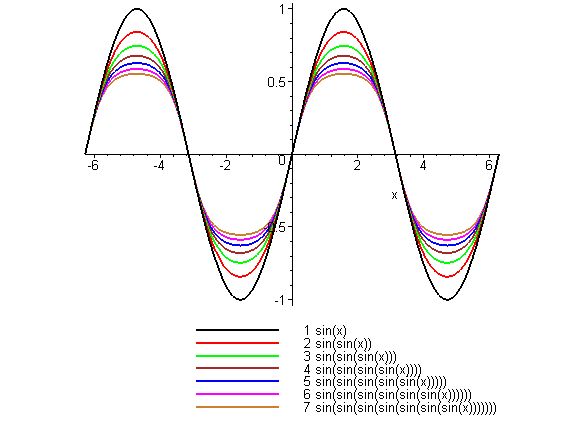}

The alternating distances between the intersection points on the cosine curve, iteral of the first order, are $l_1 = 2D$ and $L_1 = 2(\pi - D)$. The alternating distances between the intersection points on the higher order $n > 1$ iteral curves are $l_n = 2D$ and $L_n = \pi -2D$, where $\frac{\pi}{2} > 2D$. Here, $l_1 = l_n$ and $L_1 = 2L_n + l_n$. The line $y = D$ is never the center of the iteral of cosine ranges. This is the cause of the alternation.

\section{The Series of Iterated Cosine and Sine}

Using D'Alembert's ratio test it is easy to establish that both cosine and sine Maclaurin series
\begin{displaymath}
\cos(z) = 1 - \frac{z^2}{2!} + \frac{z^4}{4!} - \cdots = \sum_{n=0}^{\infty}(-1)^n\frac{z^{2n}}{(2n)!}
\end{displaymath}
\begin{displaymath}
\sin(z) = z - \frac{z^3}{3!} + \frac{z^5}{5!} - \cdots = \sum_{n=0}^{\infty}(-1)^n\frac{z^{2n+1}}{(2n+1)!}
\end{displaymath}
are \textit{absolutely convergent} for any complex argument $z$. Thus, the series representing the second power of cosine and sine expressed as the Cauchy products also converge absolutely
\begin{displaymath}
\cos^2(z) = \sum_{n=0}^{\infty}(-1)^n\frac{z^{2n}}{(2n)!}\sum_{k=0}^n
\left(
\begin{array}{c}
2n\\
2k
\end{array}
\right), \; \textrm{where} \; \left(
\begin{array}{c}
n\\
k
\end{array}
\right) = \frac{n!}{k!(n - k)!}
\end{displaymath}
\begin{displaymath}
\sin^2(z) = \sum_{n=0}^{\infty}(-1)^n\frac{z^{2n+1}}{(2n+1)!}\sum_{k=0}^n
\left(
\begin{array}{c}
2(n + 1)\\
2k+1
\end{array}
\right).
\end{displaymath}
Here, $\cos^2(z)$ and $\sin^2(z)$ are the ordinary second powers of cosine and sine. Due to the new iteral notation \cite{salov}, they cannot be misinterpreted in mixed expressions as the twice iterated trigonometric functions.

The Cauchy products sequentially created for higher powers of each function from the absolutely convergent components are absolutely convergent series on each step \cite{korn}. Our interest is even powers of cosine and odd powers of sine.

One way to get the series for $\textrm{\LARGE{\CYRI}}_{z}^{2}\cos(z) = \cos(\cos(z))$ is to replace $z$ in each power term of the Maclaurin series expressing cosine with the series itself. For instance, the first two terms are
\begin{displaymath}
\textrm{\LARGE{\CYRI}}_{z}^{2}\cos(z) = \cos(\cos(z)) = \textrm{\LARGE{\CYRI}}_{z}^{2}\sum_{n=0}^{\infty}(-1)^n\frac{z^{2n}}{(2n)!} = 1 - \frac{\cos^2(z)}{2!} + \cdots = 
\end{displaymath}
\begin{displaymath}
= 1 - \frac{\sum_{n=0}^{\infty}(-1)^n\frac{z^{2n}}{(2n)!}\sum_{k=0}^n
\left(
\begin{array}{c}
2n\\
2k
\end{array}
\right)}{2!} + \cdots = 1 - (\frac{1}{2!} - \frac{z^2}{2!} + \frac{1}{3}\frac{z^4}{2!}\cdots) + \cdots
\end{displaymath}
Due to the Cauchy product construction, for any power of cosine there is one term in the series not containing $z$ and equal to one. Combining these units divided by the corresponding factorials gives for the second order iteral of cosine the infinite series $1 - \frac{1}{2!} + \frac{1}{4!} - \frac{1}{6!} + \frac{1}{8!} - \cdots$ followed by the terms containing the even powers of $z$. Absolute convergence ensures that regrouping the terms retains the limit value. The initial series independent on $z$ is equal to $\cos(1)$.

All powers of $z$ coming from the original series for cosine will stay in the Cauchy products. Presence of units prevents their vanishing. For instance, the second powers of $z$ found in the increasing powers of $\cos^{2n}(z)$ come for $n = 1, 2, 3...$ as $-z^2, -2z^2, -3z^2, \; \ldots \; -nz^2$. After replacing $z$ in the series for $\cos$ with the same series, each of these terms is multiplied by its own factor $-\frac{1}{2!},\frac{1}{4!}, -\frac{1}{6!}, \; \ldots, \; (-1)^{n} \frac{1}{2n!}$ yielding
\begin{displaymath}
\textrm{\LARGE{\CYRI}}_{z}^{2}\cos(z) = \cos(\cos(z)) = \cos(1) + (\frac{1}{2!}-\frac{2}{4!}+\frac{3}{6!} - \cdots)z^2 + \cdots \; .
\end{displaymath}
Multiplying and dividing the terms within the intermediate brackets by two gives $\frac{1}{2}(\frac{2}{2!}-\frac{4}{4!} +\frac{6}{6!} - \cdots + (-1)^{n+1}\frac{2n}{(2n!)} - \cdots) = \frac{1}{2}\sin(1)$ and
\begin{displaymath}
\textrm{\LARGE{\CYRI}}_{z}^{2}\cos(z) = \cos(\cos(z)) = \cos(1) + \frac{1}{2}\sin(1)z^2 + higher \; even \; powers \; of \; z .
\end{displaymath}
Let us compare this series with the series for the iteral of cosine of the order one 
\begin{displaymath}
\textrm{\LARGE{\CYRI}}_{z}^{1}\cos(z) = \cos(z) = 1 - \frac{1}{2}z^2 + higher \; even \; powers \; of \; z .
\end{displaymath}
We notice that a) the first term independent on $z$ in the iteral of the second order is the first iterate of $\cos$ of the first term of the iteral of the first order; b) the new coefficient at the power $z^2$ is less by absolute value than the corresponding coefficient in the iteral of the first order; c) the sign of this coefficient alternates. The successive series may not get odd powers of $z$. While with the increasing number of iterations the coefficients containing powers of $z$ are less obvious at the moment, the term not containing $z$ is getting clear
\begin{displaymath}
\textrm{\LARGE{\CYRI}}_{z}^{n}\cos(z) = \textrm{\LARGE{\CYRI}}_{1}^{n-1}\cos(z) + terms \; containing \; z.
\end{displaymath}

Another way to the same conclusion is representing an iteral of cosine of the order $n$ by the Maclaurin series (compare it with the previous series)
\begin{displaymath}
\textrm{\LARGE{\CYRI}}_{z}^{n}\cos(z) = \textrm{\LARGE{\CYRI}}_{0}^{n}\cos(z) + \frac{d\textrm{\LARGE{\CYRI}}_{z}^{n}\cos(z)}{dz}(z_0=0)\frac{z}{1!}+\frac{d^2\textrm{\LARGE{\CYRI}}_{z}^{n}\cos(z)}{dz^2}(z_0=0)\frac{z^2}{2!}+...
\end{displaymath}
\begin{displaymath}
= \textrm{\LARGE{\CYRI}}_{1}^{n-1}\cos(z) + (-1)^n \prod_{k=0}^{n-1} \sin(\textrm{\LARGE{\CYRI}}_{0}^{k}\cos(z))\frac{z}{1!}+\frac{d^2\textrm{\LARGE{\CYRI}}_{z}^{n}\cos(z)}{dz^2}(z_0=0)\frac{z^2}{2!}+... \; ,
\end{displaymath}
where the first term is $f(z_0 = 0) = \textrm{\LARGE{\CYRI}}_{0}^{n}\cos(z) = \textrm{\LARGE{\CYRI}}_{1}^{n-1}\cos(z)$. Reduction of the order by one takes place because of the first iteration $\cos(z_0=0) = 1$. Based on the Miller's notice and our (Tumanov and author) proof
\begin{displaymath}
\lim_{n \rightarrow \infty} \textrm{\LARGE{\CYRI}}_{1}^{n-1}\cos(z) = \textrm{\LARGE{\CYRI}}_{1}^{\infty}\cos(z) = D.
\end{displaymath}
Since the limit of the iterals of cosine as it was proved does not depend on the initial \textit{real} value, the limit of the remaining terms exists and is equal to zero
\begin{displaymath}
\textrm{\LARGE{\CYRI}}_{z}^{\infty}\cos(z) = \textrm{\LARGE{\CYRI}}_{1}^{\infty}\cos(z) + \lim_{n \rightarrow \infty} terms \; containing \; z = D + 0.
\end{displaymath}
The "calculator experiment" confirms it for the real initial values. Computer experiments show that for some initial complex $z$ the iterates of cosine approach the Dottie number. However, for some other the real and imaginary parts grow and cause the overflow condition. While the expression for the first derivative of the iterals of cosine and sine as well as the series constructed in this section remain valid for the complex arguments, the limit on the iteral order does not have to hold because \textit{successive} iterals and the first derivatives are unbound for some complex arguments. To clarify details we need expressions for higher order derivatives of the iterals of cosine. A simplifying condition is that all derivatives for Maclaurin series are evaluated at $z_0 = 0$.

\section{Higher Order Derivatives}

The first derivative of the iteral of $\cos$ of the order $m > 0$ can be written as
\begin{displaymath}
\frac{d\textrm{\LARGE{\CYRI}}_{z}^{m}\cos(z)}{dz} = (-1)^m \prod_{k=0}^{m-1} \sin(\textrm{\LARGE{\CYRI}}_{z}^{k}\cos(z)) = (-1)^m \sin(z) \prod_{k=1}^{m-1} \sin(\textrm{\LARGE{\CYRI}}_{z}^{k}\cos(z))
\end{displaymath}
It is equal to zero for $z = 0$ because of the factor $\sin(z=0)$. Thus, the term $\frac{f'(z_0=0)}{1!}z$ vanishes. In the limit case $m \rightarrow \infty$ we reuse the property that the infinite product is equal to zero, if at least one of its factors is equal to zero.

The second derivative is taken as the first derivative of the product of sines composing cosine iterals of increasing orders. In courses of Calculus, it is proved using mathematical induction that $[uvw...s]' = u'vw ... s + uv'w ... s + uvw' ... s + uvw ... s'$ \cite[p. 201]{fihtengoltz}. The author writes the formulas involving the number of factoring functions $f_i(x), i = 1, ..., m$. If $f_i'(x)$ exists for all $i$, then
\begin{displaymath}
\left(\prod_{i=1}^{m}f_i(x)\right)'=\sum_{i=1}^{m} \left\{f_i'(x)\left(\prod_{j=1}^{i-1}f_j(x)\right)\left (\prod_{k=i+1}^{m}f_k(x) \right ) \right \}.
\end{displaymath}
If we put $\frac{f_i(x)}{f_i(x)} = 1$ also when $f_i(x) = 0$ and take care grouping and computing this ratio first, then simpler formulas are
\begin{displaymath}
\left(\prod_{i=1}^{m}f_i(x)\right)'=\sum_{i=1}^{m}\left\{f_i'(x)\frac{\prod_{j=1}^{m}f_j(x)}{f_i(x)}\right\}=\sum_{i=1}^{m}\left\{\frac{f_i'(x)}{f_i(x)}\prod_{j=1}^{m}f_j(x)\right \}.
\end{displaymath}
Each formula in the right side "skips" the factor $f_i(x)$ in the $i$th summand. We notice that in the last formula $\frac{f_i'(x)}{f_i(x)}=\ln(f_i(x))'$. Applying the last substitution implies that $f_i(x) \ne 0$. All summands but the first one have the factor $\sin(z)$. They are equal to zero at $z_0 = 0$. The non-vanishing term is
\begin{displaymath}
\frac{d^2\textrm{\LARGE{\CYRI}}_{z}^{m}\cos(z)}{dz^2}(z_0=0) = (-1)^m \cos(z_0=0) \prod_{k=1}^{m-1} \sin(\textrm{\LARGE{\CYRI}}_{0}^{k}\cos(z)) =
\end{displaymath}
\begin{displaymath}
= (-1)^m \prod_{k=1}^{m-1} \sin(\textrm{\LARGE{\CYRI}}_{0}^{k}\cos(z))=(-1)^m \prod_{k=1}^{m-1} \sin(\textrm{\LARGE{\CYRI}}_{1}^{k-1}\cos(z)).
\end{displaymath}
For $m=2$ we get $\sin(1)$. This completes the coefficient for $z^2$ in the Maclaurin series of the second order iteral of cosine as $\frac{1}{2}\sin(1)$ previously obtained by a different method. Let us notice that
\begin{displaymath}
\lim_{m \rightarrow \infty}\left( (-1)^m \prod_{k=1}^{m-1} \sin(\textrm{\LARGE{\CYRI}}_{1}^{k-1}\cos(z))\right)\frac{z^2}{2!} = 0.
\end{displaymath}
This is because the iteral of cosine with the initial value one tends to $D$ with $m \rightarrow \infty$ (Miller's notice). Due to the inequalities $0.6 < \sin(D) < D < 0.8$, the product forms an infinitesimal variant multiplied by $\frac{z^2}{2!}$ selected in advance and playing the role of a constant.

In order to move forward we need a formula for $n$th derivative of the product of $m$ functions, where each function has $n$th derivative. For $m > 0$  and $n=1$ the formulas are constructed. For $m=2$ and $n \geq 0$ this is the Leibniz formula
\begin{displaymath}
(u(x)v(x))^{(n)}=\sum_{i=0}^n
\left(
\begin{array}{c}
n\\
i
\end{array}
\right)
u^{(n-i)}(x)v^{(i)}(x).
\end{displaymath}
Here, the "powers" in brackets denote derivatives of the corresponding orders. The order zero means the function itself. Using binomial, trinomial, ..., $n$-nomial trees corresponding to the second, third, ..., $n$th derivative it is easy to show that the sum representing the $n$th derivative of the product of $m$ functions will contain $m^n$ summands. Each is a product of $m$ factors. Each factor is either a derivative of a certain order or a function. If the function is viewed as a derivative of the order zero, then the sum of derivative orders for each product must be equal to $n$. Some of the products will be identical and can be combined getting the multinomial contributing coefficients. This will reduce the number of summands. However, the sum of the coefficients still must be equal to $m^n$. Then, for $i=1,2,\dots,m$ and $j_i=0,1,\dots,n$
\begin{displaymath}
\left(
\prod_{i=1}^m f_i(x)\right)^{(n)} = \sum_{j_1, ..., j_m; \sum_{i=1}^m j_i = n} \left( \frac{n!}{\prod_{k=1}^m j_k!}\prod_{l=1}^m f_l^{(j_l)}(x) \right).
\end{displaymath}
While each index $j_i$ can get any value from $0$ to $n$, summation is done on all those combinations of the indexes, where the sum of their values remains equal to $n$. Let us verify that the formula works in a few already known to us cases.

If $m = 2$ and $n=0$, then it reduces to the product of two functions. In this case the only summation combination is $(j_1=0, j_2 = 0)$ and the only coefficient is equal to $\frac{0!}{0!0!}=1$. Then, $(f_1(x)f_2(x))^{(0)} = f_1(x)f_2(x)$ from the left side is equal to $1 f_1^{(0)}(x)f_2^{(0)}(x) = f_1(x)f_2(x)$ from the right side.

If $m=2$ and $n \geq 0$, then $j_2 = n - j_1$ and the relevant $n + 1$ combinations are $(j_1, n - j_1)$. The coefficients are equal to $\frac{n!}{j_1!j_2!} = \frac{n!}{j_1!(n - j_1)!} = $ binomial coefficient $(n, j_1)$. The product of functions contains only two factors $f_1^{(j_1)}(x)f_2^{(n-j_1)}$. Thus, the right side reduces to the Leibniz formula.

If $m>2$ and $n = 1$, then the relevant combinations are those, where one index is equal to one and all other are equal to zero. The total number of combinations in this case is equal to $m$. Thus, all coefficients are equal to one. Since we included the derivative of the order zero into considerations, the right side is the sum of $\prod_{l=1}^m f_l^{(j_l)}(x)$, where only one factor is the first derivative. This is exactly our result obtained earlier. The case $m = 1$ and $n > 0$ is also covered by the formula. The formula exposes similarity with expansion of $(x_1 + ... + x_m)^n$ as it was pointed to (without giving the formula) in \cite[p. 238]{fihtengoltz}. Here, $n$ without brackets is the ordinary power. The author views the formula as well known and omits the double mathematical induction proof on $m$ and $n$.

This formula should be applied to the first derivative of cosine iterals. Due to the presence of the factor $\sin(z)$ in the first derivative, each summand in the just proved formula will get the factor $\sin(z_0 = 0) = 0$ for each coefficient in front of an odd power $i$ of $\frac{z^i}{i!}$. Thus, all those terms vanish from the Maclaurin series of the iterals of cosine. This is understood, because the "last" $\cos$ series contains only even powers of $z$ and all replacements also contain only even powers of $z$. The coefficient of the Maclaurin series at $\frac{z^{2p}}{(2p)!}$ for $p > 1$ is
\begin{displaymath}
\left((-1)^m \prod_{k=0}^{m-1} \sin(\textrm{\LARGE{\CYRI}}_{z}^{k}\cos(z))\right)^{(2p-1)}=(-1)^m\left(\prod_{k=1}^{m} \sin(\textrm{\LARGE{\CYRI}}_{z}^{k-1}\cos(z))\right)^{(2p-1)}.
\end{displaymath}
The power is odd $2p-1$ because we start from the first derivative given as a product of functions. The formula for the $n$th derivative of the product gives
\begin{displaymath}
(-1)^m\left(\prod_{k=1}^{m} \sin(\textrm{\LARGE{\CYRI}}_{z}^{k-1}\cos(z))\right)^{(2p-1)}
\end{displaymath}
\begin{displaymath}
= (-1)^m \sum_{j_1, ..., j_{m}; \sum_{i=1}^{m} j_i = 2p-1} \left( \frac{(2p-1)!}{\prod_{k=1}^{m} j_k!}\prod_{l=1}^{m} \left(\sin(\textrm{\LARGE{\CYRI}}_{z}^{l-1}\cos(z))\right)^{(j_l)} \right).
\end{displaymath}
Let us prove that this coefficient tends to zero in the limit $m \rightarrow \infty$ for all $p > 1$ at $z = 0$. For $p = 1$ we have proved it already.

The first function under the product in the left side is $\sin(z)$. Its odd derivative $2p-1$ is equal to $\cos(z)$ with a precision of sign. Thus, the summand for the combination $(j_1 = 2p -1, j_2=j_3=...=j_m=0)$ remains with the first factor $\cos(z_0=0)=1$ and the multinomial coefficient equal to one independently on $m$. All other factors for $l > 1$ are sines of iterals of cosine $\sin(\textrm{\LARGE{\CYRI}}_{0}^{l-1}\cos(z))$. This construction as we have seen tends to zero with $m \rightarrow \infty$  as an infinite product of the numbers closer and closer approaching $0.6 < \sin(\textrm{D}) < 0.8$.

If for $l>1$ at least one index $j_l=1$, then the sine of the iteral of cosine being a compound function produces the factor $\sin(z)$, which will make the entire summand equal to zero at $z_0=0$. Summands with $j_1=\textrm{even number}$ will be equal to zero because of the derivative of the first factor $\sin(z)$. Summands with $j_1=0$ will be equal to zero because they will contain the factor $\sin(z)$ as the first function. We see that the question about equality of a summand to zero at $z_0=0$ translates to the question in the number theory and consideration of individual $j_i \geq 0$ within the constraining sum $j_1+...+j_i+...+j_m=2p-1$. Here, the right side is fixed for a given $p$, while the left sum gets increasing number of non-negative terms for $m \rightarrow \infty$. The infinite "majority" of them must be equal to zero leaving the function-factors "untouched". The product of the untouched factors will tend to zero at $z_0=0$ with $m \rightarrow \infty$. Without losing generality (since other cases are already considered), we set $j_1=2q-1, 1 \leq q < p$ and get $j_2 + ... +j_i + ... +j_m = 2(p - q)$. The number of "touched" (differentiated) factors is finite for any given $p$. Each of them is bound (they are computed at real $z_0 = 0$). Thus, the product of a finite number of bound factors and infinite number of factors tending to zero, should tend to zero. Presence of multinomial coefficients does not influence on this conclusion. Thus, applying the Maclaurin series we get for the terms containing $z$ in the limit of the iteral order approaching infinity the sum of the \textit{infinite} number of \textit{infinitesimal} variants. We know that it vanishes for real $z$ due to our and other proofs such as graphical analysis or Newton's method. For complex $z$ it is not so.

\section{Filled in Julia Sets of Cosine and Sine}

Using the standard set building notation and iterals, a \textit{filled in Julia set} can be defined as $z \in \mathrm{C}, \; J_c = \{v \in \mathrm{C} : \; $|\CYRI$_{v}^{\infty}f(z)| < \infty\}$ \cite{salov}. The definition of the \textit{Mandelbrot set} is another example of usefulness of the iterals: $z \in \mathrm{C}, \; M = \{c \in \mathrm{C} : \; $|\CYRI$_{0}^{\infty}(z^2 + c)| < \infty\}$ \cite{salov}. For cosine and sine the filled in Julia sets can be defined as
\begin{displaymath}
z \in \mathrm{C}, \; J_{\cos} = \{v \in \mathrm{C} : \left|\textrm{\LARGE{\CYRI}}_{v}^{\infty}\cos(z)\right| < \infty\}
\end{displaymath}
\begin{displaymath}
z \in \mathrm{C}, \; J_{\sin} = \{v \in \mathrm{C} : \left|\textrm{\LARGE{\CYRI}}_{v}^{\infty}\sin(z)\right| < \infty\}
\end{displaymath}
Computer experiments show that for some complex numbers iterals tend to infinity by absolute value. The structure of the sets is sophisticated. The points on a complex plane defining the sets form queer patterns. Appendix presents a C++ program and figure of a subset.

\paragraph{Acknowledgments.}  I would like to thank Alexander Tumanov for reviewing the article, preventing penetration of an error into it, and two ideas completing two proofs marked in the text.

\section{Appendix}
The author has written the C++ program iteral\_julia\_set.cpp. The Standard C++ Library has the class complex and a set of overloaded functions using the complex numbers \cite{stroustrup}. The comment is self-explanatory.
\begin{verbatim}
// Program scans a rectangular area of the complex plane specified
// by two points x1 + i * y1 and x2 + i * y2. The grid is the number
// of points along x- and y-sides of the rectangle. The last parameter
// selects iteral of cos or sin. Parameters should be given in the
// command line. Example:
//      iteral_julia_set -2.5 -2.5 2.5 2.5 500 cos > f.txt
// Program outputs coordinates of complex points, where an iteral tends
// to inifinity. The output format is friendly to Gnuplot program. In
// the example the output is redirected to the file f.txt. This file
// can be used with gnuplot or wgnuplot applying a single command
//      plot "f.txt" with dots 8
#include <stdexcept>    /* Standard C++ Library, exceptions */
#include <iostream>     /* Standard C++ Library, I/O streams */
#include <iomanip>      /* Standard C++ Library, manipulators */
#include <string>       /* Standard C++ Library, strings */
#include <complex>      /* Standard C++ Library, complex numbers */
#include <cmath>        /* Standard C++ Library, cos and sin */
using namespace std;
// The Iteral template for a function of one variable.
// n is the number of iterations >= 0
// v is the initial value
// f is the function such as cos, sin, etc.
template<class T, class F>
T iteral(unsigned int n, const T& v, F f)
{
    T   tmp = v;
    for(unsigned int i = 0; i < n; ++i) tmp = f(tmp);
    return tmp;
}
int main(int argc, char*argv[])
{
    try {
        if(argc < 7) {
            cout << "Usage: " << argv[0] << " x1 y1 x2 y2 grid cos|sin"
                 << endl;
            return  0;
        }
        const unsigned int n = atoi(argv[5]), ITERATIONS = 50;
        if(n < 2)
            throw   invalid_argument("Grid (" + string(argv[5])
                        + ") must be >= 2");
        if(string(argv[6]) != "sin" && string(argv[6]) != "cos")
            throw   invalid_argument("Type sin or cos but not "
                        + string(argv[6]));
        complex<double> (*f)(const complex<double>&) =
            string(argv[6]) == "cos"
            ? (complex<double> (*)(const complex<double>&))cos
            : (complex<double> (*)(const complex<double>&))sin;
        const complex<double> z1(atof(argv[1]), atof(argv[2]));
        const complex<double> z2(atof(argv[3]), atof(argv[4]));
        const double dzr = (z2 - z1).real() / (n - 1);
        const complex<double> dzi(0.0, (z2 - z1).imag() / (n - 1));
        complex<double> z(z1);
        for(unsigned int r = 0; r < n; ++r, z += dzr) {
            for(unsigned int i = 0; i < n; ++i, z += dzi) {
                if(norm(iteral(ITERATIONS, z, f)) < 10)
                    cout << setw(25) << setprecision(16) << z.real()
                            << " " << setw(25) << setprecision(16)
                            << z.imag() << endl;
            }
            z = complex<double>(z.real(), z1.imag());
        }
    }
    catch(const exception& e) {
        cerr    << e.what() << endl;
        return  -1;
    }
    catch(...) {
        cerr    << "Unknown exception" << endl;
        return  -2;
    }
    return  0;
}
\end{verbatim}
For each point in a rectangular grid the program evaluates iterals of cosine or sine for 50 iterations. If the iteral value is less than 10 by absolute value, then corresponding coordinates of the point are output. Otherwise, the point is skipped in the output. The standard output can be redirected to a text file and plotted using gnuplot or wgnuplot shareware. The result of plotting $500 * 500 = 250000$ scanned and checked points in the square $-2.5-2.5i$ and $2.5+2.5i$ are on figure. Setting ranges and grid of the rectangle works as zoom useful for investigation of fine details of the iterals.

\includegraphics[width=130mm]{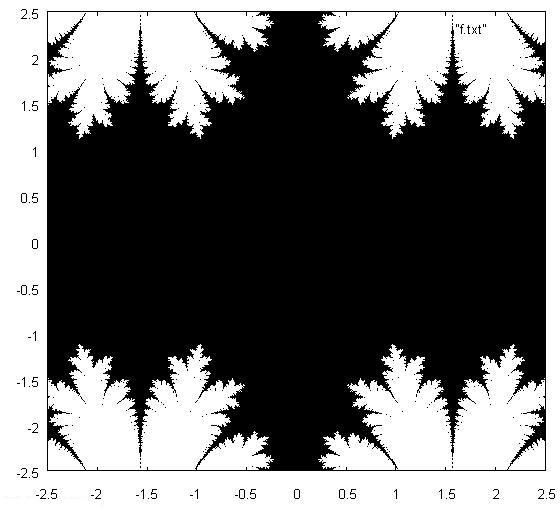}

\bigskip

\noindent\textbf{Valerii Salov} received his M.S. from the Moscow State University, Department of Chemistry in 1982 and his Ph.D. from the Academy of Sciences of the USSR, Vernadski Institute of Geochemistry and Analytical Chemistry in 1987.  He is the author of the articles on analytical, computational, and physical chemistry, the book Modeling Maximum Trading Profits with C++, \textit{John Wiley and Sons, Inc., Hoboken, New Jersey}, 2007, the paper on iterals of functions \cite{salov}, and papers in \textit{Futures Magazine}.

\noindent\textit{v7f5a7@comcast.net}

\end{document}